\documentclass[12pt]{amsart}

\newtheorem{THM}{Theorem}[section]
\newtheorem{LMA}[THM]{Lemma}
\newtheorem{PROP}[THM]{Proposition}
\newtheorem{CORO}[THM]{Corollary}

\newtheorem{POR}[THM]{Porism}

\numberwithin{equation}{section}

\usepackage{amsmath}
\usepackage{amssymb}
\usepackage{euscript}

\newcommand{\showon}{\begin{eqnarray*}}
\newcommand{\showoff}{\end{eqnarray*}}

\newcommand{\drop}{\smallsetminus}

\newcommand{\goesto}{\rightarrow}
\newcommand{\one}{\boldsymbol{1}}
\newcommand{\zero}{\boldsymbol{0}}

 \renewcommand{\a}{\mathbf{a}}
\renewcommand{\b}{\mathbf{b}}
 \newcommand{\CC}{\mathbb{C}}
\renewcommand{\c}{\mathbf{c}}

\renewcommand{\L}{\EuScript{L}}
\newcommand{\M}{\EuScript{M}} 
\newcommand{\N}{\EuScript{N}}

\newcommand{\U}{\EuScript{U}}\renewcommand{\u}{\mathbf{u}}
\renewcommand{\v}{\mathbf{v}}
\newcommand{\W}{\EuScript{W}}\newcommand{\w}{\mathbf{w}}

\newcommand{\Y}{\EuScript{Y}} \newcommand{\y}{\mathbf{y}}

\begin{document}

\title[Matroid inequalities]{Remarks on one combinatorial application
of the Aleksandrov--Fenchel inequalities}

\author{David G. Wagner}
\address{Department of Combinatorics and Optimization\\
University of Waterloo\\
Waterloo, Ontario, Canada\ \ N2L 3G1}
\email{\texttt{dgwagner@math.uwaterloo.ca}}


\begin{abstract}
In 1981, Stanley applied the Aleksandrov--Fenchel Inequalities to prove a
logarithmic concavity theorem for regular matroids.  Using ideas from electrical
network theory we prove a generalization of this for the wider class of
matroids with the ``half--plane property''.  Then we explore a nest of inequalities
for weighted basis--generating polynomials that are related to these ideas.
As a first result from this investigation we find that every matroid of rank
three or corank three satisfies a condition only slightly weaker than the 
conclusion of Stanley's theorem.
\end{abstract}

\maketitle

\begin{center}
\emph{Dedicated with great admiration\\
to Richard Stanley\\
on the occasion of his 60th birthday.}
\end{center}

\section{Introduction.}

In 1981, Stanley \cite{St} applied the Aleksandrov--Fenchel Inequalities
to prove the following logarithmic concavity result:
\begin{THM}[Stanley, 1981]
Let $\M$ be a matroid on a set $E$.  Let 
$\pi=(S,T,C_{1},\ldots,C_{k})$ be an ordered partition of $E$ into
pairwise disjoint nonempty subsets, and fix nonnegative integers
$c_{1},\ldots,c_{k}$.  For each $0\leq j\leq |S|$, let
$M_{j}(\pi)$ be the number of bases $B$ of $\M$ such that $|B\cap S|=j$
and $|B\cap C_{i}|=c_{i}$ for all $1\leq i\leq k$. 
If $\M$ is regular then for each $1\leq j\leq |S|-1$,
$$\frac{M_{j}(\pi)^{2}}{\binom{|S|}{j}^{2}}\geq
\frac{M_{j-1}(\pi)}{\binom{|S|}{j-1}}\cdot
\frac{M_{j+1}(\pi)}{\binom{|S|}{j+1}}.$$
\end{THM}
\noindent
Stanley's proof proceeds by constructing a set of zonotopes with the
desired mixed volumes.

A few years later, answering a question raised by Stanley,
Godsil \cite{Go} strengthened this conclusion as follows:
\begin{THM}[Godsil, 1984]
With the hypotheses and notation of Theorem 1.1, the polynomial
$\sum_{j=0}^{|S|} M_j(\pi) x^j$ has only real (nonpositive) zeros.
\end{THM}
\noindent
The well--known Newton's Inequalities (item (51) of \cite{HLP}) show that
Theorem 1.2 implies Theorem 1.1.  Godsil's proof employs a determinantal
theorem used by Schneider \cite{Sch} to prove the Aleksandrov--Fenchel Inequalities.

The aim of this paper is to publicize a recent extension of these results
(Theorem 4.5 of \cite{CW}), to collect the scattered details of its proof into a
self--contained whole, and to present some preliminary findings
on related inequalities.  It can also be regarded as an idiosyncratic
introduction to the 
literature on the half--plane property, Rayleigh monotonicity, and related\
concepts \cite{BB,C2,COSW,CW,FB,FM,Ly,Pe,SW,W}.  We extend Theorem {1.2} in two
directions -- by relaxing the hypothesis
and by strengthening the conclusion.  For the hypothesis, we replace the condition that
the matroid is regular by the weaker condition that the matroid has the
half--plane property explained in Section 2.  We strengthen the
conclusion by introducing positive real weights $\y:=\{y_e:\ e\in E\}$ on the
elements of the ground--set $E(\M)$.  The weight of a basis $B$ of $\M$ is then
$\y^B:=\prod_{e\in B} y_e$, and with the notation of Theorem 1.1 we let
$M_j(\pi,\y):=\sum_B \y^B$ with the sum over all bases $B$ of $\M$ such that
$|B\cap S|=j$ and $|B\cap C_{i}|=c_{i}$ for all $1\leq i\leq k$. 
\begin{THM}[Theorem 4.5 of \cite{CW}]
With the above hypotheses and notation, if $\M$ has the half--plane property
then the polynomial $$\sum_{j=0}^{|S|} M_j(\pi,\y) x^j$$ has only real (nonpositive)
zeros.
\end{THM}

\section{Electrical Networks and Matroids.}

Our proof of Theorem 1.3 builds on ideas originating with Kirchhoff's formula
for the effective conductance of a (linear, resistive, direct current) electrical network.
Such a network is represented by a graph $G=(V,E)$ with real positive weights
$\y:=\{y_e:\ e\in E\}$ specifying the conductance of each edge of the graph.
We use the notation
$$G(\y):=\sum_T \y^T$$
for the sum of $\y^T:=\prod_{e\in T} y_e$ over all spanning trees $T$ of $G$,
and $\y>\zero$ to indicate that every edge--weight is positive.
Fixing two vertices $v,w\in V$, the effective conductance of the network measured
between $v$ and $w$ is -- by Kirchhoff's Formula \cite{Ki} --
$$\Y_{vw}(G;\y)=\frac{H^f(\y)}{H_f(\y)}.$$
Here $H$ denotes the graph obtained from $G$ by adjoining a new edge $f$ with
ends $v$ and $w$, $H^f=G$ is $H$ with $f$ deleted, and $H_f$ is $H$ with $f$
contracted.

One intuitive property of electrical networks is that if the conductance $y_e$
of some edge $e$ of $G$ is increased, then $\Y_{vw}(G;\y)$ does not decrease.
That is,
$$\frac{\partial}{\partial y_e}\Y_{vw}(G;\y)\geq 0.$$
This property is known as \emph{Rayleigh monotonicity}.
After some calculation using $H_f(\y)=H_f^e(\y)+y_e H_{ef}(\y)$ \emph{et cetera},
this is found to be equivalent to the condition that if
$\y>\zero$ then
$$H_e^f(\y)H_f^e(\y)\geq H_{ef}(\y)H^{ef}(\y).$$
(The deletion/contraction notation is extended in the obvoius way.)

A less obvious property of the effective conductance is that if every $y_e$ is
a complex number with positive real part then  $\Y_{vw}(G;\y)$ is a complex
number with nonnegative real part.  Physically, this corresponds to the fact that
if every edge of an alternating current circuit dissipates energy then the whole
network cannot produce energy.
Some minor hijinx with the equation $H(\y)=H^f(\y)+y_f H_f(\y)$ shows that
this is equivalent to the condition that
if Re$(y_e)>0$ for all $e\in E(H)$ then $H(\y)\neq 0$. Such a polynomial $H(\y)$
is said to have the \emph{half--plane property}.

The combinatorics of the preceeding three paragraphs carries over \emph{mutatis
mutandis} to matroids in general.  In place of a graph $G$ we have a matroid $\M$.
The edge--weights $\y$ become weights on the ground--set $E(\M)$ of $\M$.
In place of the spanning tree generating function $G(\y)$ we have the basis
generating function
$$M(\y):=\sum_{B\in\M} \y^B.$$
Since this is insensitive to loops we might as well think of $\M$ as given by its
set of bases.  For disjoint subsets $I,J\subseteq E$ the contraction of $I$ and
deletion of $J$ from $\M$ is given by
$$\M_I^J:=\{B\drop I:\ B\in\M\ \mathrm{and}\  I\subseteq B\subseteq E\drop J\}.$$

Rayleigh monotonicity corresponds to the inequalities
$$M_e^f(\y)M_f^e(\y)\geq M_{ef}(\y)M^{ef}(\y)$$
for all $e,f\in E$ and $\y>\zero$.  A matroid satisfying these inequalities is
called a \emph{Rayleigh matroid}.  If the basis generating polynomial $M(\y)$
has the half--plane property then we say that $\M$ has the \emph{half--plane property},
or is a \emph{HPP matroid}.

\section{The Half--Plane Property.}

Our first item of business is to show that every regular matroid is a
HPP matroid.  For graphs, Proposition $3.1$
is part of the ``folklore'' of electrical engineering.  We take it
from Corollary 8.2(a) and Theorem 8.9 of \cite{COSW}, but include
the short and interesting proof for completeness.

A matrix $A$ of complex numbers is a \emph{sixth--root of unity}
matrix provided that every nonzero minor of $A$ is a sixth--root
of unity.  A matroid $\M$ is a \emph{sixth--root of unity matroid}
provided that it can be represented over the complex numbers by a
sixth--root of unity matrix.  For example, every regular matroid
is a sixth--root of unity matroid.  Whittle \cite{Wh} has shown that
a matroid is a sixth--root of unity matroid if and only if it is
representable over both $GF(3)$ and $GF(4)$. (It is worth noting that
Godsil's proof can be adapted to prove Theorem 1.3 in the special case of
sixth--root of unity matroids.)

\begin{PROP}
Every sixth--root of unity matroid is a HPP matroid.
\end{PROP}
\begin{proof}
Let $A$ be a sixth--root of unity matrix of full row--rank $r$, 
representing the matroid $\M$, and let $A^{*}$ denote the 
conjugate transpose of $A$.  Index the columns of $A$ by
the set $E$, and let $Y:=\mathrm{diag}(y_{e}:\ e\in E)$ be a
diagonal matrix of indeterminates.  For an $r$--element subset $S
\subseteq E$, let $A[S]$ denote the square submatrix of $A$ supported on
the set $S$ of columns.  By the Binet--Cauchy formula,
$$\det(AYA^{*})=\sum_{S\subseteq E:\ |S|=r}|\det A[S]|^{2}\y^{S}=M(\y)$$
is the basis--generating polynomial of $\M$,  since $|\det A[S]|^{2}$
is $1$ or $0$ according to whether or not $S$ is a basis of $\M$.

Now we claim that if Re$(y_{e})>0$ for all $e\in E$, then $AYA^{*}$ is
nonsingular.  This suffices to prove the result.  Consider any nonzero
vector $\v\in\CC^{r}$.  Then $A^{*}\v\neq\zero$ since the columns of
$A^{*}$ are linearly independent.  Therefore
$$\v^{*}AYA^{*}\v=\sum_{e\in E}y_{e}|(A^{*}\v)_{e}|^{2}$$
has strictly positive real part, since for all $e\in E$ the numbers
$|(A^{*}\v)_{e}|^{2}$ are nonnegative reals and at least one of these
is positive.  In particular, for any nonzero $\v\in\CC^{r}$, the
vector $AYA^{*}\v$ is nonzero.  It follows that $AYA^{*}$ is 
nonsingular, completing the proof.
\end{proof}
The same proof shows that for any complex matrix $A$ of full row--rank
$r$, the polynomial
$$\det(AYA^{*})=\sum_{S\subseteq E:\ |S|=r}|\det A[S]|^{2}\y^{S}$$
has the half--plane property.  The weighted analogue of Rayleigh 
monotonicity in this case is discussed from a probabilistic point of 
view by Lyons \cite{Ly}.  It is a surprising fact that a complex
matrix $A$ of full row--rank $r$ has $|\det A[S]|^{2}=1$ for all
nonzero rank $r$ minors if and only if $A$ represents a sixth--root of
unity matroid (Theorem $8.9$ of \cite{COSW}).

The proof of Theorem 1.3 is based on the following lemmas regarding the
half--plane property.  The paper \cite{COSW} gives a much more thorough
development of the theory.

\begin{LMA}
Let $P(\y)$ be a polynomial in the variables
$\y=\{y_{e}:\ e\in E\}$, let $e\in E$, and let the degree of $y_e$ in
$P$ be $n$.  If $P(\y)$ has the half--plane property then
$y_e^n P(\{y_f:\ f\neq e\},1/y_e)$ has the half--plane property.
\end{LMA}
\begin{proof}
This follows immediately from the fact that Re$(1/y_e)>0$ if and only if
Re$(y_e)>0$.
\end{proof}

\begin{LMA}[\cite{COSW}, Proposition 2.8]
Let $P(\y)$ be a polynomial in the variables
$\y=\{y_{e}:\ e\in E\}$.  For any $e\in E$, if $P$ has the half--plane
property then $\partial P/\partial y_e$ has the half--plane property.
\end{LMA}
\begin{proof}
Fix complex values with positive real parts for every $y_f$ with $f\in E\drop\{e\}$.
The result of these substitutions is a univariate polynomial $F(y_e)$ all the
roots of which have nonpositive real part.  Thus
$$F(y_e)=C\prod_{j=1}^n(y_e+\theta_j)$$
with Re$(\theta_j)\geq 0$ for all $1\leq i\leq n$.  It follows that if Re$(y_e)>0$ then
the real part of
$$\frac{F'(y_e)}{F(y_e)}=\sum_{j=1}^n\frac{1}{y_e+\theta_j}$$
is also strictly positive.  In particular $F'(y_e)\neq 0$.  It follows
that $\partial P/\partial y_e$ has the half--plane property.
\end{proof}

\begin{CORO}[\cite{FB}, Theorem $18$, or \cite{COSW}, Proposition $3.4$.]
Let $P(\y)$ be a polynomial in the variables
$\y=\{y_{e}:\ e\in E\}$, fix $e\in E$, and let
$P(\y)=\sum_{j=0}^{n}P_{j}(\{y_f:\ f\neq e\})y_{e}^{j}$.
If $P$ has the half--plane property then each $P_{j}$ has the 
half--plane property.
\end{CORO}
\begin{proof}
Let $n$ be the degree of $P(\y)$ in the variable $y_e$.
Let $A:=\partial^j P/\partial y_e^j$,
$B:=y_e^{n-j}A(\{y_f:\ f\neq e\},1/y_e)$, and
$C:=\partial^{n-j-1} B/\partial y_e^{n-j-1}$.
Then  $C(\y)$ is a nonzero multiple of $P_j$, and has the half--plane
property by Lemmas 3.2 and 3.3.
\end{proof}

\begin{LMA}[\cite{COSW}, Proposition $5.2$]
Let $P(\y)$ be a homogeneous polynomial in the variables
$\y=\{y_{e}:\ e\in E\}$.  For sets of nonnegative real numbers
$\a=\{a_{e}:\ e\in E\}$ and $\b=\{b_{e}:\ e\in E\}$, let
$P(\a x+\b)$ be the polyomial obtained by substituting
$y_{e}=a_{e}x+b_{e}$ for each $e\in E$.  The following are
equivalent:\\
\textup{(a)} $P(\y)$ has the half--plane property;\\
\textup{(b)} for all sets of nonnegative real numbers $\a$ and $\b$,
$P(\a x + \b)$ has only real (nonpositive) zeros.
\end{LMA}
\begin{proof}
To see that (a) implies (b), suppose that $\xi$ is a zero of $P(\a x+\b)$
that is not a nonpositive real.  Then there are complex numbers $z$ and $w$ with
positive real part such that $z/w=\xi$.  If $P(\y)$ is homogeneous of degree $r$ then
$P(\a z+\b w)= w^r P(\a\xi +\b)=0$, showing that $P(\y)$ fails to have the
half--plane property.

To see that (b) implies (a), consider any set of values $\{y_e:\ e\in E\}$
with Re$(y_e)>0$ for all $e\in E$.  There are complex numbers $z$ and $w$
with positive real parts such that all the $y_e$ are in the convex cone generated by
$z$ and $w$.  That is, for each $e\in E$ there are nonnegative reals $a_e$ and $b_e$
such that $y_e=a_e z+b_e w$.  Now
$P(\y)=w^r P(\a\xi+\b)$ in which $\xi=z/w$ is not in the interval $(-\infty,0]$,
and so $P(\y)\neq 0$.
\end{proof}

We can now prove Theorem 1.3.
\begin{proof}[Proof of Theorem 1.3.]
Let $\M$ be a HPP matroid and fix $\y>\zero$.
Let $s$, $t$, and $z_{1},\ldots,z_{k}$ be 
indeterminates, and for $e\in E$ put
$$u_{e}:=\left\{\begin{array}{ll}
y_{e}s & \mathrm{if}\ \  e\in S,\\
y_{e}t & \mathrm{if}\ \  e\in T,\\
y_{e}z_{i} & \mathrm{if}\ \ e\in C_{i}.\end{array}\right.$$
Then $M(\u)$ is a homogeneous polynomial with the half--plane property
in the variables $s,t,z_{1},\ldots,z_{k}$.  By repeated application
of Corollary 3.4, the coefficient $M_{\c}(s,t)$ of
$z_{1}^{c_{1}}\cdots z_{k}^{c_{k}}$ in
$M(\u)$ also has the half--plane property, and is homogeneous.  In fact,
$$M_{\c}(s,t) = \sum_{j=0}^{|S|}M_{j}(\pi,\y)s^{j}t^{d-j},$$
in which $d=\mathrm{rank}(\M)-(c_{1}+\cdots+c_{k})$.
Upon substituting $s=x$ and $t=1$ in $M_{\c}(s,t)$, Lemma 3.5
implies that $\sum_{j=0}^{|S|}M_{j}(\pi,\y)x^{j}$ has only real
(nonpositive) zeros, as claimed.
\end{proof}

\section{Between HPP and Rayleigh.}

This updated version of Stanley's theorem provides a link between the half--plane
property and Rayleigh monotonicity in the context of matroids.
Only the $k=0$ case of the theorem is needed.
(In fact the $k=0$ case is equivalent to the general case by various properties
of HPP matroids.)  For a subset $S\subseteq E(\M)$ and natural number $j$, let
$M_j(S,\y):=\sum_{B}\y^{B}$, with the sum over all bases $B$ of $\M$
such that $|B\cap S|=j$.  For each positive integer $m$, consider the
following conditions on a matroid $\M$:\\

\noindent
\textbf{RZ[$m$]:}\ \ If $\y>\zero$ then for all
$S\subseteq E$ with $|S|\leq m$ the polynomial $\sum_{j=0}^{|S|}
M_j(S,\y)x^{j}$ has only real zeros.\\

\noindent
\textbf{BLC[$m$]:}\ \  If $\y>\zero$ then for all
$S\subseteq E$ with $|S|=n\leq m$ and all $1\leq j\leq n-1$,
$$M_j(S,\y)^2 \geq \left[1+\frac{n+1}{j(n-j)}\right] M_{j-1}(S,\y) M_{j+1}(S,\y).$$
\vspace{2mm}

\noindent
The mnemonics are for ``real zeros'' and ``binomial logarithmic concavity'',
respectively.  We also say that a matroid satisfies RZ if it satisfies
RZ[$m$] for all $m$, and that it satisfies BLC if it satisfies
BLC[$m$] for all $m$.  Our BLC is a weighted strengthening of Stanley's ``Property P''.

The $k=0$ case of Theorem 1.3 implies that a HPP matroid satisfies RZ, and
Newton's Inequalities show that RZ[$m$] implies BLC[$m$] for every $m$.
The implications RZ[$m$]\
$\Longrightarrow$ RZ[$m-1$] and BLC[$m$]\ $\Longrightarrow$ BLC[$m-1$]
are trivial, as are the conditions RZ[$1$] and BLC[$1$].  Thus, the
weakest nontrivial condition among these is BLC[$2$].  This is in fact
\textbf{equivalent} to Rayleigh monotonicity, as the remarks after Theorem 4.3 show.

\begin{PROP}[\cite{CW}, Corollary 4.9]
Every HPP matroid is a Rayleigh matroid.
\end{PROP}
\begin{proof}
Theorem 1.3 shows that every HPP matroid satisfies BLC[2].  We show here
that if $\M$ satisfies BLC[2] then it is Rayleigh.  So, let $\y>\zero$
be positive weights on $E(\M)$ and let $S=\{e,f\}\subseteq E$.  To
prove the Rayleigh inequality $M_e^f(\y)M_f^e(\y)\geq M_{ef}(\y)M^{ef}(\y)$
it suffices to consider the case in which both $e$ and $f$ are neither
loops nor coloops.  In this case, define
another set of weights by $w_e:=M_f^e(\y)$ and $w_f:= M_e^f(\y)$ and
$w_g:=y_g$ for all $g\in E\drop\{e,f\}$.  Then, since $\w>0$ and $\M$
satisfies BLC[2], the inequality
$$M_1(S,\w)^2\geq 4 M_0(S,\w) M_2(S,\w)$$
holds.  This can be expanded to
$$(w_e M_e^f(\w)+w_f M_f^e(\w))^2 \geq 4 M^{ef}(\w)w_e w_f M_{ef}(\w),$$
and finally to
$$4(M_e^f(\y) M_f^e(\y))^2 \geq 4 M_{ef}(\y)M^{ef}(\y)M_e^f(\y)M_f^e(\y).$$
Cancellation of common (positive) factors from both sides yields the
desired inequality.
\end{proof}

In view of the implication BLC[2] $\Longrightarrow$ Rayleigh, it is interesting to
look for conditions (other than HPP) which imply BLC[$m$] for various $m$.
This is further motivated by Stanley's application of ``Property P'' to Mason's
Conjecture -- see Theorem 2.9 of \cite{St}.
The following hierarchies of \emph{strict root--binomial logarithmic concavity}
and \emph{strict logarithmic concavity} conditions are also interesting:\\

\noindent
$\sqrt{\mathbf{BLC}}$\textbf[$m$]\ \  If $\y>\zero$ then for all
$S\subseteq E$ with $|S|=n\leq m$ and all $1\leq j\leq n-1$,
if $M_j(S,\y)\neq 0$ then
$$M_j(S,\y)^2 > \left[1+\frac{1}{\min(j,n-j)}\right]
M_{j-1}(S,\y) M_{j+1}(S,\y).$$
\vspace{2mm}

\noindent
\textbf{SLC[$m$]:}\ \  If $\y>\zero$ then for all
$S\subseteq E$ with $|S|\leq m$ and all $1\leq j\leq |S|-1$,
if $M_j(S,\y)\neq 0$ then
$$M_j(S,\y)^2 > M_{j-1}(S,\y) M_{j+1}(S,\y).$$
\vspace{2mm}

We also say that a matroid satisfies
$\sqrt{\mathrm{BLC}}$ if it satisfies $\sqrt{\mathrm{BLC}}[m]$ for all $m$,
and satisfies SLC if it satisfies SLC[$m$] for all $m$.
The inequalities
$$1+\frac{1}{\min(j,n-j)}<1+\frac{n+1}{j(n-j)}\leq
\left[1+\frac{1}{\min(j,n-j)}\right]^2$$
for $1\leq j\leq n-1$ show that BLC[$m$] implies $\sqrt{\mathrm{BLC}}[m]$
for every $m$, and motivate this somewhat odd terminology.	Clearly
$\sqrt{\mathrm{BLC}}[m]$ implies SLC[$m$] for every $m$. 

In the following calculations we will usually omit explicit reference to the
variables $\y$ unless a particular substitution must be emphasized.
For $\M$ a matroid, $S\subseteq E(\M)$, and $k$ a natural number, let
$$\Psi_k M S:=\sum_{A\subseteq S:\ |A|=k} M_A^{S\drop A}M_{S\drop A}^A.$$
For example,
$$\Psi_2 M\{a,b,c\}:=M_{ab}^c M_c^{ab}+M_{ac}^b M_b^{ac}+M_{bc}^a M_a^{bc}.$$
Notice that in general if $|S|=n$ then $\Psi_k M S= \Psi_{n-k} MS$.
The Rayleigh inequality is that
$$\Psi_1 M\{e,f\}\geq 2\Psi_2 M\{e,f\}.$$
This suggests several possible generalizations, among which we will concentrate here
on the following.  For each integer $k\geq 1$ and real $\lambda>0$, say that
$\M$ is \emph{$k$--th level Rayleigh of strength $\lambda$} provided that:\\

\noindent
$\lambda$\textbf{-Ray[$k$]:}\ \ If $\y>\zero$ then for every $S\subseteq E(\M)$
with $|S|=2k$,
$$\Psi_k M S \geq \lambda\Psi_{k+1} M S.$$

The condition $2$-Ray[$1$] is exactly the Rayleigh condition.
In general, each term on the LHS of $\lambda$-Ray[$k$] occurs
twice.  Proposition 4.6 below shows that $(1+1/k)$-Ray[$k$] is an especially
natural strength for these conditions.  Interestingly, this lies right between
two of the most \textbf{useful} strengths for these conditions.  As an example, the
inequality for $(3/2)$-Ray[$2$] is
\showon
& &4\left[ M_{ab}^{cd} M_{cd}^{ab}+M_{ac}^{bd} M_{bd}^{ac}
+M_{ad}^{bc} M_{bc}^{ad}\right]\\
&\geq& 3\left[ M_{bcd}^a M_a^{bcd} + M_{acd}^b M_b^{acd}
+  M_{abd}^c M_c^{abd} + M_{abc}^d M_d^{abc} \right].
\showoff

\begin{LMA}
\textup{(a)}\
For each $k\geq 1$ and $\lambda>0$, the class of matroids satisfying
$\lambda$-$\mathrm{Ray}[k]$ is closed by
taking duals and minors.\\
\textup{(b)}\
For each $m\geq 1$, the class of matroids satisfying $\mathrm{BLC}[m]$ is closed by
taking duals and minors.\\
\textup{(c)}\
For each $m\geq 1$, the class of matroids satisfying $\sqrt{\mathrm{BLC}}[m]$ is closed by
taking duals and minors.\\
\textup{(d)}\
For each $m\geq 1$, the class of matroids satisfying $\mathrm{SLC}[m]$ is closed by
taking duals and minors.\\
\end{LMA}
\begin{proof}[Sketch of proof.]
For the matroid $\M^*$ dual to $\M$ we have $M^*(\y)=\y^E M(\one/\y)$.
From this it follows that for $0\leq k\leq n=|S|$,
$$\Psi_k M^* S(\y) = \left(\y^{E\drop S}\right)^2 \Psi_{n-k} M S(\one/\y).$$
Since $\y>\zero$ is arbitrary, one sees that $\lambda$-Ray[$k$] for $\M$ implies
$\lambda$-Ray[$k$] for $\M^*$.  Similarly
$$M_j(S,\y)^2>M_{j-1}(S,\y)M_{j+1}(S,\y)$$
implies
$$(\y^{E})^2 M^*_{n-j}(S,\one/\y)^2>
(\y^{E})^2 M^*_{n-j+1}(S,\one/\y)M^*_{n-j-1}(S,\one/\y).$$
From this it follows that SLC[$m$] for $\M$ implies
SLC[$m$] for $\M^*$.  Analogous arguments work for BLC$[m]$ and $\sqrt{\mathrm{BLC}}[m]$.

For a set $S\subseteq E(\M)$ and $g\in E\drop S$, we have
$$\Psi_k M S(\y) = y_g^2 \Psi_k M_g S(\y) + y_g Q(\y) + \Psi_k M^g S(\y)$$
for some polynomial $Q(\y)$.
If $\M$ satisfies $\lambda$-Ray[$k$] then taking the limit as $y_g\goesto 0$ shows that
the deletion $\M^g$ satisfies $\lambda$-Ray[$k$].  Similarly, multiplying by $1/y_g^2$
and taking the limit as $y_g\goesto\infty$ shows that the contraction
$\M_g$ satisfies $\lambda$-Ray[$k$].  The case of a general minor is obtained by iteration
of these two cases.   Analogous arguments work for BLC$[m]$,  $\sqrt{\mathrm{BLC}}[m]$,
and SLC$[m]$.
\end{proof}

\begin{THM} If $\M$ satisfies $2$-$\mathrm{Ray}[1]$ and 
$(1+1/k)^2$-$\mathrm{Ray}[k]$ for all $2\leq k\leq m$ then $\M$ satisfies
$\mathrm{BLC}[2m+1]$.
\end{THM}
\begin{proof}
The proof uses the following elementary inequality:\  for $N\geq 2$
real numbers $R_1,...,R_N$,
\begin{eqnarray}
(R_1+\cdots+R_N)^2\geq \frac{2N}{N-1}\sum_{\{i,j\}\subseteq\{1,...,N\}} R_i R_j.
\end{eqnarray}

Assume that $\M$ satisfies the hypothesis, and fix positive real weights
$\y>\zero$.  To show that $\M$ satisfies BLC[$2m+1$], consider a subset
$S\subseteq E$ with $|S|=n\leq 2m+1$ and an index $1\leq j\leq n-1$.
We must show that
$$M_j(S)^2 \geq \left[1+\frac{n+1}{j(n-j)}\right] M_{j-1}(S)M_{j+1}(S).$$
From the definition we have
$$M_j(S)=\sum_{A\subset S:\ |A|=j} \y^A M_A^{S\drop A}.$$
The inequality (4.1) implies that
\begin{eqnarray}
M_j(S)^2\geq \frac{2N}{N-1}\sum_{\{A,B\}} \y^A \y^B M_A^{S\drop A} M_B^{S\drop B},
\end{eqnarray}
with $N:=\binom{n}{j}$ and the sum over all pairs of distinct $j$--element subsets
of $S$.  We collect terms on the right side according to the intersection $I:=
A\cap B$ and union $U:=A\cup B$ of the indexing pair of sets $\{A,B\}$.  With 
$i:=|I|$ and $k:=j-i$ we have $|U\drop I|=2k$ and the sum of the terms on the RHS of
(4.2) with this fixed $I$ and $U$ is
\showon
& &
\frac{N}{N-1} \y^I\y^U \sum_{C\subset U\drop I:\ |C|=k}
\left(M_I^{S\drop U}\right)_{C}^{U\drop I\drop C}
\left(M_I^{S\drop U}\right)_{U\drop I\drop C}^{C}\\
&=& \frac{N}{N-1}\y^I\y^U\Psi_k M_I^{S\drop U}(U\drop I).
\showoff
Therefore
$$M_j(S)^2\geq
\frac{N}{N-1}\sum_{k=1}^{h}\sum_{(I,U)} \y^I\y^U \Psi_k M_I^{S\drop U}(U\drop I)$$
in which $h:=\min(j,n-j)$ and the inner sum is over all pairs of sets
$I\subset U\subseteq S$ with $|U|+|I|=2j$ and $|U\drop I|=2k$.
Let $\lambda_1:=2$ and $\lambda_k:=(1+1/k)^2$ for all $k\geq 2$.
Since $k\leq h\leq m$ we are assuming that $\M$ satisfies $\lambda_k$-Ray[$k$], and hence
each minor $\M_I^{S\drop U}$ satisfies $\lambda_k$-Ray[$k$] by Lemma 4.2(a).
It follows that
\begin{eqnarray}
M_j(S)^2 &\geq&
\frac{N}{N-1}\sum_{k=1}^{h}\lambda_k\sum_{(I,U)}
\y^I\y^U \Psi_{k+1}M_I^{S\drop U}(U\drop I).
\end{eqnarray}

If $h\geq 3$ then for all $1\leq k\leq h$,
$$\lambda_k\geq \lambda_h=\left(1+\frac{1}{h}\right)^2\geq 1+\frac{n+1}{j(n-j)}.$$
If $h=2$ and $n\geq 5$ then
$$\lambda_2=\frac{9}{4}>2=\lambda_1\geq  1+\frac{n+1}{j(n-j)}.$$
In these cases we conclude from (4.3) that
\showon
M_j(S)^2 &\geq&
\frac{N}{N-1}\left[1+\frac{n+1}{j(n-j)}\right]
\sum_{k=1}^{h}\sum_{(I,U)}\y^I\y^U \Psi_{k+1}M_I^{S\drop U}(U\drop I)\\
&=& \frac{N}{N-1}\left[1+\frac{n+1}{j(n-j)}\right]M_{j-1}(S)M_{j+1}(S).
\showoff
This implies the desired inequality in these cases.

If $h=1$ then either $j=1$ or $j=n-1$, so $k=1$ and $N=n$, and we conclude from
(4.3) that
\showon
M_j(S)^2 &\geq&
\frac{2n}{n-1}\sum_{(I,U)}\y^I\y^U \Psi_{2}M_I^{S\drop U}(U\drop I)\\
&=& \frac{2n}{n-1}M_{j-1}(S)M_{j+1}(S).
\showoff
This is the desired inequality when $j=1$ or $j=n-1$.

If $n\leq 3$ then $h=1$, so the only remaining case is $n=4$ and $h=2$.
In this case $j=2$ and $N=6$, and since $\lambda_2=9/4>2=\lambda_1$, from
(4.3) we conclude that
\showon
M_2(S)^2 &\geq&
\frac{6}{5}
\sum_{k=1}^{2}2\sum_{(I,U)}\y^I\y^U \Psi_{k+1}M_I^{S\drop U}(U\drop I)\\
&\geq& \frac{9}{4}M_{1}(S)M_{3}(S).
\showoff
This is the desired inequality in this case.

This completes the verification that $\M$ satisfies BLC[$2m+1$].
\end{proof}
Notice that Proposition 4.1 and the $m=1$ case of Theorem 4.3 show that the conditions 
BLC$[3]$, BLC$[2]$, and $2$-Ray$[1]$ are in fact equivalent -- this is part of Theorem
4.8 of \cite{CW}.

By examining the proof of Theorem 4.3 one sees the following.
\begin{POR}\textup{(a)}\ If $\M$ satisfies $(1+1/k)$-$\mathrm{Ray}[k]$ for
all $1\leq k\leq m$ then $\M$ satisfies $\sqrt{\mathrm{BLC}}[2m+1]$.\\
\textup{(b)}\  If $\M$ satisfies $1$-$\mathrm{Ray}[k]$ for all $1\leq k\leq m$
then $\M$ satisfies $\mathrm{SLC}[2m+1]$.
\end{POR}

\begin{CORO}
Let $\M$ be a matroid of rank $r$ with $|E|=\ell$ elements.\\
\textup{(a)}\  If $\M$ satisfies $2$-$\mathrm{Ray}[1]$ and $(1+1/k)^2$-$\mathrm{Ray}[k]$
for all $2\leq k\leq r-1$ then $\M$ satisfies \textup{BLC}.\\
\textup{(b)}\  If $\M$ satisfies $(1+1/k)$-$\mathrm{Ray}[k]$
for all $1\leq k\leq r-1$ then $\M$ satisfies $\sqrt{\mathrm{BLC}}$.\\
\textup{(c)}\
If $\M$ satisfies $1$-$\mathrm{Ray}[k]$ for all $1\leq k\leq \min(r-1,(\ell-2)/4)$
then $\M$ satisfies \textup{SLC}.\\
\end{CORO}
\begin{proof}
If $\M$ has rank $r$ then $\M$ satisfies any condition $\lambda$-Ray$[k]$ trivially
for all $k\geq r$ because for any $S\subseteq E(\M)$ with $|S|=2k$, $\Psi_{k+1} M S=0$
identically.  In part (a) Theorem 4.3 thus implies that $\M$ satisfies
BLC[$2m+1$] for all $m$ -- that is BLC.  In part (b) Porism 4.4(a) supports the
analogous implication.  Similarly, in part (c) if
$r-1\leq(\ell-2)/4$ then Porism 4.4(b) implies that $\M$ satisfies SLC.
If $(\ell-2)/4<r-1$ then Porism 4.4(b) implies
that $\M$ satisfies SLC[$\lfloor \ell/2\rfloor$].  For any $S\subset E$ with $|S|=n\leq
\ell/2$ and any $1\leq j\leq n$ we have $M_j(S)=M_{r-j}(E\drop S)$.  The SLC inequalities
for the set $S$ thus imply the SLC inequalities for the set $E\drop S$.  It follows that
$\M$ satisfies SLC.
\end{proof}

One of course wants examples of matroids satisfying these higher level
Rayleigh conditions.  So far there are no substantial results in this direction.
Direct computations (with the aid of \textsc{Maple}) and \emph{ad hoc}
arguments have shown that the Kuratowski graphs $\mathsf{K}_5$ and
$\mathsf{K}_{3,3}$ satisfy $(3/2)$-Ray$[2]$ but not $(9/4)$-Ray$[2]$.
The hypothesis of Theorem 4.3 thus seems to be extremely restrictive.
The corresponding calculation for $\mathsf{K}_6$ is so far intractible.
Similar computations for  $\lambda$-Ray$[k]$ with $k\geq 3$ have not yet been tried.

One easy sufficient condition is the following naively optimistic generalization of the
Rayleigh condition.
\begin{PROP}
Assume that $\M$ is such that for all disjoint subsets $A$ and $B$ of
$E(\M)$ with $|A|=|B|=k$ and $b\in B$, if $\y>\zero$ then
$$\M_A^B M_B^A \geq  M_{Ab}^{B\drop b} M_{B\drop b}^{Ab}.$$
Then $\M$ satisfies $(1+1/k)$-$\mathrm{Ray}[k]$.
\end{PROP}
\begin{proof}
Let $\M$ be as in the hypothesis, let $\y>\zero$, and
fix a set $S\subseteq E(\M)$ with $|S|=2k$.  To verify the inequality
$k\Psi_k MS\geq (k+1)\Psi_{k+1} MS$, form a graph $G$ with bipartition $(X,Y)$
as follows.  The vertices in $X$ are the $k$--element subsets of $S$ and the
vertices of $Y$ are the $(k+1)$ element subsets of $S$.  There is an edge from
$A\in X$ to $A'\in Y$ whenever $A\subset A'$.  Thus, every vertex of $X$ has degree
$k$ and every vertex of $Y$ has degree $k+1$.  To each edge $\{A,Ab\}$ of $G$
associate the weight $M_A^{S\drop A}M_{S\drop A}^A$.  The total weight assigned to
the edges is thus
$$k\sum_{A\in X}M_A^{S\drop A}M_{S\drop A}^A=k\Psi_k MS.$$
On the other hand, for each $A'\in Y$ the sum of the weights of edges incident with
$A'$ is
\showon
& &\sum_{b\in A'} M_{A'\drop b}^{S\drop (A'\drop b)}M_{S\drop (A'\drop b)}^{A'\drop b}\\
&\geq& \sum_{b\in A'} M_{A'}^{S\drop A'} M_{S\drop A'}^{A'}
=(k+1)M_{A'}^{S\drop A'} M_{S\drop A'}^{A'}.
\showoff
It follows that the sum of the weights of the edges of $G$ is at least
$(k+1)\Psi_{k+1} MS$, completing the proof.
\end{proof}

Unfortunately, it is not hard to find \textbf{planar graphs} which fail
to satisfy the hypothesis of Proposition 4.6 in the case $k=2$.
An easy example is the $4$--wheel $\mathsf{W}_4$
formed by a cycle of four ``rim'' edges with the vertices joined to a new ``hub'' vertex
by four ``spoke'' edges.  The spoke edges are labelled $a$, $b$, $c$, $d$ in
cyclic order, and the rim edges have weights $y$, $1$, $y$, $1$ in cyclic
order starting with the edge forming a triangle with $a$ and $b$.
With $\W$ denoting the graphic matroid of
$\mathsf{W}_4$, we have $W_{ab}^{cd}=W_{cd}^{ab}=2y+1$ and $W_{abc}^d=y+1$
and $W_d^{abc}=2y(y+1)$, so that the inequality
$$W_{ab}^{cd}W_{cd}^{ab}\geq W_{abc}^d W_d^{abc}$$ is equivalent to
$$(2y+1)^2\geq 2y(y+1)^2.$$
This inequality $2y^3\leq 2y+1$ is not satisfied for all $y>0$, so that
$\W$ does not meet the hypothesis of Proposition 4.4.  Note that $\W$ does satisfy
$(3/2)$-Ray$[2]$, however, since $\mathsf{W}_4$ is a minor of $\mathsf{K}_5$.

In general, the condition $(1+1/k)$-Ray[$k$] asserts that the inequalities in 
the hypothesis of Proposition 4.6 hold ``locally on average'' in some 
sense.  

To close this section we show that the condition $\sqrt{\mathrm{BLC}}$ has
consequences for Mason's Conjecture \cite{Ma} -- this is inspired by Theorem
2.9 of Stanley \cite{St}.  For a matroid $\M$ of rank $r$ and $0\leq j\leq r$,
let $I_j(\M)$ denote the number of $j$--element independent sets of $\M$.  Mason's
Conjecture is that for any matroid, $I_j^2\geq I_{j-1}I_{j+1}$ for all
$1\leq j\leq r-1$.  Let $\M\oplus\U_{\ell,\ell}$ denote the free extension of $\M$
by $\ell$ points, and let $T_r(\M\oplus\U_{\ell,\ell})$ denote the truncation of 
this matroid down to rank $r$.

\begin{THM}
Let $\M$ be a matroid of rank $r$.  Assume that, for infinitely many integers
$\ell\geq 1$, $T_r(\M\oplus\U_{\ell,\ell})$ satisfies $\sqrt{\mathrm{BLC}}[\ell]$.
Then $\M$ satisfies Mason's Conjecture.
\end{THM}
\begin{proof}
Let $\L:=T_r(\M\oplus\U_{\ell,\ell})$ for some $\ell\geq 1$, and let $S:=E(\L)\drop
E(\M)$.  For any $0\leq j\leq r$ we have $L_j(S,\one)=\binom{\ell}{j}I_{r-j}(\M)$.
If $1\leq j\leq r-1$ and $\ell\geq 2j$ is such that $\L$ satisfies
$\sqrt{\mathrm{BLC}}[\ell]$ then
$$ L_j(S,\one)^2 > \left[1+\frac{1}{j}\right]L_{j-1}(S,\one)L_{j+1}(S,\one).$$
From this we obtain
$$I_{r-j}(\M)^2 > \frac{\ell-j}{\ell-j+1}I_{r-j+1}(\M)I_{r-j-1}(\M).$$
Since this holds for a sequence of $\ell\goesto\infty$ we conclude that
$$I_{r-j}(\M)^2 \geq I_{r-j+1}(\M)I_{r-j-1}(\M).$$
As this holds for all $1\leq j\leq r-1$, $\M$ satisfies Mason's Conjecture.
\end{proof}

\section{Matroids of rank three.}

In \cite{W} it is shown that every matroid of rank (at most) three satisfies the
Rayleigh condition $2$-Ray$[1]$.

\begin{THM}
Every matroid of rank three satisfies $(3/2)$-$\mathrm{Ray}[2]$.
\end{THM}
\begin{proof}
Let $\M$ be a matroid of rank three and let $\y>\zero$ be positive weights on
$E(\M)$. Consider any four--element subset $S=\{a,b,c,d\}$ of $E$.  There are
three cases (up to permuting the elements of $S$):\\
(i)\ $\{a,b,c,d\}$ are collinear;\\
(ii)\ $\{a,b,c\}$ are collinear and $d$ is not on this line;\\
(iii)\ $\{a,b,c,d\}$ are in general position.

In case (i) it is easy to see that $\Psi_3 M S=0$, so that the inequality
$2\Psi_2 MS \geq 3\Psi_3 MS$ is trivially satisfied.

For cases (ii) and (iii) we consider the coefficient of each monomial in
$\Psi_2 MS$ and in $\Psi_3 MS$ in turn.

First consider monomials of the form $y_e^2$ with $e\in E\drop S$. 
The coefficient of $y_e^2$ in $\Psi_2 MS$ is nonnegative.
On the other hand, the coefficient of $y_e^2$ in $M_{bcd}^aM_a^{bcd}$ is zero,
and similarly the coefficient of $y_e^2$ in $\Psi_3 MS$ is zero.  It follows that
in both cases (i) and (ii), the coefficient of  $y_e^2$ in $2\Psi_2 MS-3\Psi_3 MS$
is nonnegative.

\setlength{\unitlength}{0.8mm}
\begin{figure}
\begin{picture}(150,230)
\thicklines

\put(10,210){\line(1,0){40}}
\put(10,210){\circle*{2}}\put(9,204){$2$}
\put(20,210){\circle*{2}}\put(19,204){$3$}
\put(30,210){\circle*{2}}\put(29,204){$4$}
\put(40,210){\circle*{2}}\put(39,204){$5$}
\put(50,210){\circle*{2}}\put(49,204){$6$}
\put(30,225){\circle*{2}}\put(33,224){$1$}
\put(5,225){I.}

\put(100,205){\line(1,0){45}}
\put(100,205){\line(3,2){30}}
\put(100,205){\circle*{2}}\put(99,199){$6$}
\put(115,205){\circle*{2}}\put(114,199){$5$}
\put(130,205){\circle*{2}}\put(129,199){$4$}
\put(145,205){\circle*{2}}\put(144,199){$2$}
\put(115,215){\circle*{2}}\put(111,217){$3$}
\put(130,225){\circle*{2}}\put(133,224){$1$}
\put(95,225){II.}

\put(10,165){\line(1,0){45}}
\put(10,165){\circle*{2}}\put(9,159){$3$}
\put(25,165){\circle*{2}}\put(24,159){$4$}
\put(40,165){\circle*{2}}\put(39,159){$5$}
\put(55,165){\circle*{2}}\put(54,159){$6$}
\put(20,180){\circle*{2}}\put(23,179){$1$}
\put(45,180){\circle*{2}}\put(48,179){$2$}
\put(5,185){III.}

\put(100,155){\line(1,0){40}}
\put(100,155){\line(2,3){20}}
\put(110,170){\line(2,-1){30}}
\put(120,155){\line(0,1){30}}
\put(100,155){\circle*{2}}\put(99,149){$3$}
\put(120,155){\circle*{2}}\put(119,149){$5$}
\put(140,155){\circle*{2}}\put(139,149){$2$}
\put(120,165){\circle*{2}}\put(123,166){$4$}
\put(110,170){\circle*{2}}\put(106,172){$6$}
\put(120,185){\circle*{2}}\put(123,184){$1$}
\put(95,185){IV.}

\put(10,105){\line(1,0){40}}
\put(10,105){\line(2,3){20}}
\put(30,135){\line(2,-3){20}}
\put(10,105){\circle*{2}}\put(9,99){$2$}
\put(30,105){\circle*{2}}\put(29,99){$4$}
\put(50,105){\circle*{2}}\put(49,99){$3$}
\put(20,120){\circle*{2}}\put(16,122){$6$}
\put(40,120){\circle*{2}}\put(43,121){$5$}
\put(30,135){\circle*{2}}\put(29,138){$1$}
\put(5,135){V.}

\put(100,105){\line(1,0){40}}
\put(100,105){\line(2,3){20}}
\put(100,105){\circle*{2}}\put(99,99){$6$}
\put(120,105){\circle*{2}}\put(119,99){$5$}
\put(140,105){\circle*{2}}\put(139,99){$2$}
\put(110,120){\circle*{2}}\put(106,122){$4$}
\put(135,125){\circle*{2}}\put(138,124){$3$}
\put(120,135){\circle*{2}}\put(123,134){$1$}
\put(95,135){VI.}

\put(10,80){\line(1,0){40}}
\put(10,60){\line(1,0){40}}
\put(10,80){\circle*{2}}\put(9,74){$1$}
\put(10,60){\circle*{2}} \put(9,54){$2$}
\put(30,80){\circle*{2}}\put(29,74){$3$}
\put(30,60){\circle*{2}} \put(29,54){$5$}
\put(50,80){\circle*{2}}\put(49,74){$4$}
\put(50,60){\circle*{2}} \put(49,54){$6$}
\put(5,85){VII.}

\put(100,60){\line(1,0){40}}
\put(100,75){\circle*{2}} \put(99,69){$1$}
\put(100,60){\circle*{2}} \put(99,54){$4$}
\put(120,80){\circle*{2}}\put(119,74){$2$}
\put(120,60){\circle*{2}} \put(119,54){$5$}
\put(140,75){\circle*{2}} \put(139,69){$3$}
\put(140,60){\circle*{2}} \put(139,54){$6$}
\put(95,85){VIII.}

\put(65,40){\circle*{2}}\put(68,39){$1$}
\put(85,40){\circle*{2}}\put(88,39){$2$}
\put(95,25){\circle*{2}}\put(98,24){$3$}
\put(85,10){\circle*{2}}\put(88,9){$4$}
\put(65,10){\circle*{2}}\put(68,9){$5$}
\put(55,25){\circle*{2}}\put(58,24){$6$}
\put(50,40){IX.}

\end{picture}
\caption{The six--element rank three matroids.}
\end{figure}
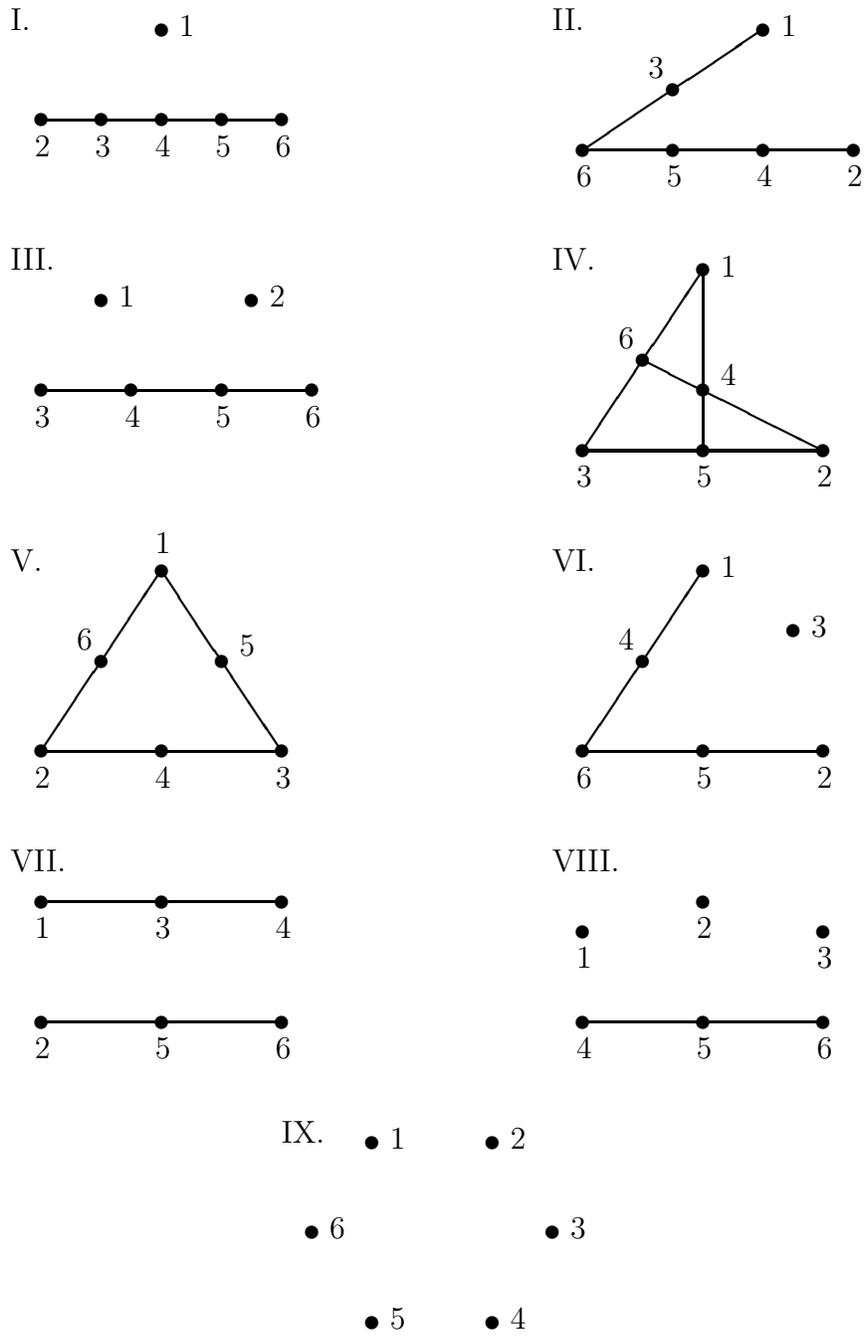

Next consider monomials of the form $y_ey_f$ with $e\neq f$ in $E\drop S$.
The coefficiens of $y_ey_f$ in $\Psi_2 MS$ and in $\Psi_3 MS$ depend only on the
restriction $\N:=\M|\{a,b,c,d,e,f\}$ and the position of $\{a,b,c,d\}$ within $\N$.
Figure $1$ lists all the rank three matroids on six elements, up to isomorphism.
Table $1$ lists the possible cases for $\N$ and $\{a,b,c,d\}$, up to isomorphism,
in case (ii).
Table $2$ lists the possible cases for $\N$ and $\{a,b,c,d\}$, up to isomorphism,
in case (iii).  In each table the first column specifies $\N$ and $\{a,b,c,d\}$,
the second column gives the coefficient of $y_ey_f$ in $\Psi_2 MS$, and the
third column gives the coefficient of $y_ey_f$ in $\Psi_3 MS$.  The fourth column computes
the coefficient of the corresponding monomial of $2\Psi_2 MS - 3\Psi_3 MS$. 
Only in case (iii) subcase $\mathrm{IV}\{1,2,3,4\}$ do we obtain a negative coefficient.

\begin{table}
$$\begin{array}{|r|r|r|c|}
\hline
\N\{a,b,c,d\} & \Psi_2 & \Psi_3 & 2\Psi_2-3\Psi_3\\ \hline

\mathrm{I}   \{2,3,4,1\} & 0  & 0 & 0 \\
\mathrm{II}  \{1,3,6,2\} & 8  & 2 & 10 \\
\mathrm{II}  \{2,4,5,1\} & 6  & 3 & 3 \\
\mathrm{II}  \{2,4,6,1\} & 6  & 3 & 3 \\
\mathrm{III} \{3,4,5,1\} & 6  & 3 & 3 \\
\mathrm{IV}  \{1,3,6,2\} & 8  & 2 & 10 \\
\mathrm{V}   \{1,2,6,3\} & 8  & 3 & 7 \\
\mathrm{V}   \{1,2,6,4\} & 10 & 3 & 11 \\
\mathrm{VI}  \{1,4,6,2\} & 10 & 3 & 11 \\
\mathrm{VI}  \{1,4,6,3\} & 12 & 3 & 18 \\
\mathrm{VII} \{1,3,4,2\} & 12 & 3 & 15 \\
\mathrm{VIII}\{4,5,6,1\} & 12 & 3 & 15 \\
\hline
\end{array}$$
\caption{Monomials of shape $y_e y_f$ in case (ii).}
\end{table}

\begin{table}
$$\begin{array}{|r|r|r|c|}
\hline
\N\{a,b,c,d\} & \Psi_2 & \Psi_3 & 2\Psi_2-3\Psi_3\\ \hline

\mathrm{II}  \{1,2,3,4\} & 8  & 2 & 10 \\
\mathrm{III} \{1,2,3,4\} & 8  & 2 & 10 \\
\mathrm{IV}  \{1,2,3,4\} & 4  & 4 & -4 \\
\mathrm{V}   \{1,2,4,5\} & 6  & 4 & 0 \\
\mathrm{V}   \{1,4,5,6\} & 8  & 3 & 7 \\
\mathrm{VI}  \{1,2,4,5\} & 8  & 4 & 4 \\
\mathrm{VI}  \{1,2,3,4\} & 10 & 3 & 11 \\
\mathrm{VI}  \{1,2,3,6\} & 10 & 4 & 8 \\
\mathrm{VII} \{1,2,3,5\} & 8  & 4 & 4 \\
\mathrm{VIII}\{1,2,4,5\} & 10 & 4 & 8 \\
\mathrm{IX}  \{1,2,3,4\} & 12 & 4 & 12\\
\hline
\end{array}$$
\caption{Monomials of shape $y_e y_f$ in case (iii).}
\end{table}

Given a set $S=\{a,b,c,d\}$ in general position in $\M$, there are at most three
points in $E\drop S$ that can appear in a restriction of $\M$ containing $S$ in subcase
$\mathrm{IV}\{1,2,3,4\}$;  these are the points $e:=\overline{ab}\cap\overline{cd}$,
$f:=\overline{ac}\cap\overline{bd}$, and $g:=\overline{ad}\cap\overline{bc}$.
Notice that (when they exist) each of $y_e^2$, $y_f^2$ and $y_g^2$ occurs in
$2\Psi_2 MS-3\Psi_3 MS$ with coefficient $8$.  If at most one of $e$, $f$, and $g$
exists in $E(\M)$ then subcase $\mathrm{IV}\{1,2,3,4\}$ does not arise.
If two exist -- say $e$ and $f$ -- then
$$8y_e^2+8y_f^2-4y_ey_f=6y_e^2+6y_f^2+2(y_e-y_f)^2$$
can be used to absorb the term with negative coefficient into a square of a binomial.
If all of $e$, $f$, and $g$ exist in $E(\M)$ then
\begin{eqnarray*}
& & 8y_e^2+8y_f^2+8y_g^2-4y_ey_f-4y_ey_g-4y_fy_g\\
&=& 4y_e^2+4y_f^2+4y_g^2+2(y_e-y_f)^2+2(y_e-y_g)^2+2(y_f-y_g)^2
\end{eqnarray*}
can be used to absorb the terms with negative coefficients into a sum of squares
of binomials.

In each of cases (i), (ii), and (iii) we can thus write $2\Psi_2 MS-3\Psi_3 MS$
as a positive sum of monomials and squares of binomials in the variables $\y$.
It follows that $\M$ satisfies the condition $(3/2)$-Ray$[2]$.
\end{proof}

Theorem 5.1, the result of \cite{W}, Lemma 4.2(c), and Corollary 4.5(b) immediately
imply the following.

\begin{THM}
Every matroid of rank three or corank three satisfies $\sqrt{\mathrm{BLC}}$.
\end{THM}
As seen in \cite{COSW}, there are many matroids of rank three that are not
HPP matroids -- the Fano and Pappus matroids are two familiar examples.
Theorem 5.2 begins to explore situations in which the hypothesis of Theorem 1.3
does not apply, but something close to its conclusion does hold.  Unfortunately
Theorem 5.2 implies nothing new about Mason's Conjecture, which is trivial in
rank three.

\end{document}